\documentclass{amsart}

\usepackage{amsmath}
\usepackage{amsfonts}
\usepackage{amsopn}
\usepackage{amssymb}

\newcommand{\R}{\ensuremath{\mathbb{R}}}
\newcommand{\C}{\ensuremath{\mathbb{C}}}

\newcommand{\Z}{\ensuremath{\mathbb{Z}}}

\DeclareMathOperator{\imag}{Im}

\DeclareMathOperator{\area}{area}

\DeclareMathOperator{\Diff}{Diff}
\DeclareMathOperator{\rk}{rk}

\def\cprime{$'$}

\newtheorem{thm}{Theorem}

\newtheorem{lemma}[thm]{Lemma}
\newtheorem{prop}[thm]{Proposition}
\newtheorem{cor}[thm]{Corollary}
\title{Scaled Relators and Dehn Functions for Nilpotent Groups}
\author{Robert Young}
\address{Department of Mathematics\\
University of Chicago\\
5734 S. University Avenue\\
Chicago, Illinois 60637\\
U.S.A.}
\date{\today}
\email{rjyoung@math.uchicago.edu}

\newcommand{\etalchar}[1]{$^{#1}$}

\providecommand{\MR}{\relax\ifhmode\unskip\space\fi MR }
\providecommand{\href}[2]{#2}

\begin{document}
\bibliographystyle{amsalpha}
\begin{abstract}  
A homogeneous nilpotent Lie group has a scaling automorphism
determined by a grading of its Lie algebra.  Many proofs of upper
bounds for the Dehn function of such a group depend on being able to
fill curves with discs compatible with this grading; the area of such
discs changes predictably under the scaling automorphism.  In this
paper, we present combinatorial methods for finding such bounds.
Using this method, we give new proofs of some results on Dehn
functions of nilpotent groups, prove theorems on central powers and
certain quotients of nilpotent groups, and construct the first example
of a torsion-free nilpotent group of class 3 with a quadratic Dehn
function.
\end{abstract}
\maketitle

\section{Introduction}
Dehn functions sit at a juncture between geometric and combinatorial
group theory; a key theorem of Gromov, proved in \cite{Bri} and \cite{BurTab},
connects the Dehn function of a group with the Dehn function of the
universal cover of a compact space with that fundamental group.  In
the case of nilpotent groups, both combinatorial and geometric methods
have proved productive in proving isoperimetric inequalities,
partially because nilpotent discrete groups and nilpotent Lie groups
are closely linked.  In particular, Mal\cprime cev's Embedding
Theorem states that torsion-free nilpotent groups are cocompact
lattices in nilpotent Lie groups, so the Dehn function of a
torsion-free nilpotent group is equivalent to the Dehn function of a
nilpotent Lie group.  

In this paper, we will develop some methods for finding isoperimetric
inequalities for lattices in homogeneous Lie groups that combine the
geometry of the Lie group and the combinatorics of a presentation of
the lattice.  Briefly, one way to find isoperimetric inequalities for
nilpotent groups, developed by Gromov in \cite{GroAI} and \cite{GroCC}, is to fill curves with discs
whose areas change predictably under the scaling automorphism.  We
show that in the case of homogeneous Lie groups, it suffices to do
this for the relators in the presentation of a lattice.
Much of the utility of this technique comes from the fact that many
nilpotent groups have presentations with relators that can be easily
filled by discs.

We will first summarize results about Dehn functions for nilpotent
groups.  In general, one can bound the Dehn function of a nilpotent
group in terms of its nilpotency class: a finitely presented nilpotent
group of class $k$ satisfies an isoperimetric inequality
$\delta(n)<n^{k+1}$.  Gromov suggested a proof involving
infinitesimally invertible operators in \cite{GroAI}; a combinatorial
proof can be found in \cite{GeHoRi}.  This bound is not sharp; for
instance, the higher-dimensional Heisenberg groups are class 2
nilpotent groups with quadratic isoperimetric inequalities.  This fact
was claimed by Thurston\cite{E} and proofs can be found in
\cite{Allcock,GroCC,SapOls}; it also follows from our result on central powers in section
\ref{appssec}.  In section \ref{appssec}, we will also present an
example of a class 3 nilpotent group with a quadratic isoperimetric
inequality.  It is still an open question how far the general bound
can be from sharp.  

The best lower bounds known for Dehn functions of nilpotent groups
come from using the centralized isoperimetric function, an analog of
signed area which is easier to compute than the isoperimetric
function; see \cite{BaMiSh} or \cite{Pittet}.  In many cases, this
bound combined with the general bound above gives a sharp estimate of
the Dehn function; one such case is for free nilpotent groups.  Few
other sources of lower bounds for nilpotent groups are known and it is
unknown whether this centralized isoperimetric function gives sharp
bounds for nilpotent groups in general.  In the case of finitely
presented groups, the bounds it gives are not sharp; in \cite{BaMiSh},
there is an example of a Baumslag-Solitar group for which the
centralized isoperimetric function is linear and the Dehn function is
exponential.

In Section \ref{preliminaries}, we will define the Dehn function of a
group and give some background on the geometry of nilpotent Lie
groups.  In Section \ref{mainthsec}, we will state
and prove our main theorem and give some corollaries, and in Section
\ref{appssec}, we will apply it to several examples, including central
powers of nilpotent groups and a nilpotent group of class 3 with
quadratic Dehn function.

The author would like to thank Tim Riley for his helpful comments on a draft of this paper and
Mark Sapir for the conversation that inspired Corollary
\ref{quotdehn}.

\section{Preliminaries}\label{preliminaries}
We begin by defining the Dehn function of a group and of a manifold; for a survey on Dehn functions, see \cite{Bri}.
Let $H=\{a_1,\dots,a_d|r_1,\dots,r_s\}$ be a group.  Any word $w$ in
the $a_i$ representing the identity can be reduced to the trivial word
by using the relators; that is, there is a decomposition
\begin{equation}\label{worddecomp}w=\prod_{i=1}^k g_i^{-1} r_{b_i}^{\pm1} g_i,\end{equation}
where the equality is
taken in the free group generated by the $a_i$.  We define the area
$\delta(w)$ of $w$ to be the smallest $k$ for which such a
decomposition exists.

Equivalently, take the universal cover of the presentation 2-complex.
This cover is a simply connected 2-complex whose 1-skeleton is the Cayley
graph of $H$ and whose faces correspond to conjugates of relators.
Words then correspond to closed curves in the complex and
decompositions as in (\ref{worddecomp}) correspond to van Kampen
diagrams(roughly, mappings of discs into the 2-complex) for $w$.  The
area $\delta(w)$ of $w$ is then the minimal number of faces in such a
van Kampen diagram.

We then define the {\em Dehn function} of $H$ by setting 
$$\delta(n)=\max\delta(w),$$ 
where the maximum is taken over all words representing the identity
with at most $n$ letters.  This function depends on the presentation
of $H$, but only by a linear distortion of the domain and range.  We define the partial
ordering 
$$f\prec g\text{ iff $\exists A,B,C,D,E$ s.t. $f(n)\le A g(B n + C)+D
n +E$ for all $n$}$$ 
and let $f\equiv g$ iff $f\prec g$ and $f\succ g$.  Then if $\delta$
and $\delta'$ are Dehn functions corresponding to two finite
presentations of the same group, $\delta\equiv\delta'$, so we can
talk about the Dehn function of a group up to this equivalence
relation.  The equivalence class of the Dehn function is in fact a quasi-isometry invariant.  That is, if $\delta$ and $\delta'$ are Dehn functions of quasi-isometric groups, then $\delta\equiv\delta'$.  

If $M$ is a simply-connected Riemannian manifold, one can define a similar
geometric Dehn function.  Given a Lipschitz closed curve $\gamma$ in
$M$, one can define its area $\delta(\gamma)$ to be the infimum of the
areas of all the Lipschitz discs bounded by $\gamma$ and thus define a
Dehn function of $M$: 
\begin{equation}\label{Dehnmandef}\delta(n)=\sup\delta(\gamma),\end{equation}
where the supremum is taken over all Lipschitz closed curves of length
at most $n$.

The connection between these two definitions is given in a theorem stated by Gromov and proved in \cite{Bri} and \cite{BurTab}:
\begin{thm}
If $M$ is a simply connected Riemannian manifold and $H$ is a finitely presented group acting properly discontinuously, cocompactly, and by isometries on $M$, then $\delta_H\equiv \delta_M$.
\end{thm}
Thus, to understand the Dehn function of a lattice in a nilpotent Lie
group, we will consider discs embedded in the Lie group.  In
particular, we will consider lattices with presentations whose Cayley
graphs embed naturally into the Lie group.  Using either
proof of the theorem, one can easily show that taking the maximum area over curves
contained in the Cayley graph in definition (\ref{Dehnmandef}) gives
an equivalent definition of the Dehn function.

We will give some preliminaries on the geometry of homogeneous
nilpotent Lie groups.  Let $G$ be a nilpotent Lie group and $\Gamma$
its Lie algebra.  Then if we let $\Gamma_1=\Gamma$ and
$\Gamma_i=[\Gamma,\Gamma_{i-1}]$, there is a filtration
$$\Gamma=\Gamma_1\supset\dots\supset\Gamma_k\supset\Gamma_{k+1}=0$$
corresponding to the lower central series, and
$[\Gamma_i,\Gamma_j]\subset \Gamma_{i+j}$.  If $\Gamma_k\ne 0$, we
call $k$ the {\em nilpotency class} of $G$.

If $G$ is simply connected and there is a grading $$\Gamma=V_1\oplus\dots\oplus V_k$$ such that
$[V_i,V_j]\subset V_{i+j}$, then we call $G$
{\em homogeneous}.  We define a scaling automorphism $s_t$ of $\Gamma$ so
that $s_t$ acts on $V_i$ by multiplication by $t^i$ and extend it to
an automorphism of $G$.  

We will use scaling automorphisms to approximate curves by simpler
curves, so it will be important to know how the area of a curve
changes under scaling.  We put a left-invariant Riemannian metric
on $G$ and extend the $V_i$ to left-invariant subbundles of $TG$.
Then $s_t$ distorts vectors in $TG$ according to where they lie in the
grading.  If $v\in V_1\oplus\dots\oplus V_i$, then $||s_{t*}(v)||<c
t^i||v||$ for some $c$.  Call a piecewise $C^2$ map whose tangent vectors lie in
$V_1\oplus\dots\oplus V_i$ {\em $(V_1\oplus\dots\oplus
V_i)$-horizontal} or {\em horizontal} when the subspace is $V_1$.
Then if $f:D^2\to G$ is $(V_1\oplus\dots\oplus V_i)$-horizontal,
then $\area(s_t\circ f)<c^2 t^{2i}\area f$.

A left-invariant Riemannian metric is not always the best choice of
metric for a homogeneous nilpotent Lie group.  In particular, this
metric is distorted by scaling automorphisms.  To avoid this problem,
we will often use a left-invariant subriemannian or
Carnot-Carath\'{e}odory metric instead.  To construct this metric,
we first call a curve horizontal if it is always tangent to $V_1$; one can
show that any two points in the group can be connected by a horizontal
curve.  By choosing a metric on $V_1$, one can measure the length of
horizontal curves, and we define the distance between two points to be
the infimum of the lengths of horizontal curves connecting them.
Since the tangent vectors to the curves are all in $V_1$, $s_t$
changes the length of the curve by a factor of $t$, so the scaling
automorphisms act as homotheties on the metric.

If $H$ is a lattice in $G$ which is generated by $a_1,\dots,a_d$ such
that $a_i=\exp(v_i)$ for $v_i\in V_1$, we say that $H$ is compatible
with the grading.  Then if $n$ is an integer, $s_n(H)$ is generated by
$a_1^n,\dots,a_d^n$ and is a subgroup of $H$.  Similarly, if $\langle
a_1,\dots,a_d|r_1,\dots,r_s\rangle$ is a presentation where
$a_i=\exp(v_i)$ for $v_i\in V_1$, we call it a presentation compatible
with the grading.  

\begin{prop}\label{compatlat}
If $G$ is a homogeneous nilpotent group with a grading 
$$\Gamma=V_1\oplus\dots\oplus V_k$$ 
of its Lie algebra and $H\subset G$ is a lattice, then $G$ has a lattice $H_0$ compatible with the grading.
\end{prop}
\begin{proof}  
Define $p_i:\Gamma\to V_i$ to be the projection onto the $i$th factor
of the grading and define $H^{(1)}=H$, $H^{(i)}=[H,H^{(i-1)}]$ to be
the lower central series of $H$.  Consider the subgroup $H_0$ of $G$
generated by $\exp(p_1(\log H))$.  If this is a lattice, it is
compatible with the grading by definition.  We begin by showing
inductively that $p_j(\log H_0^{(j)})\subset p_j(\log H^{(j)})$.
Since each of the $p_j(\log H^{(j)})$ is a lattice in $V_j$, this will
show that $H_0$ is discrete.

In the case $j=1$, note that, by the Baker-Campbell-Hausdorff formula
\begin{align*}
\exp(u) \exp(v) &= \exp(u + v + \frac{1}{2} [u, v] +  \dots)\\
p_1(\log(\exp(u) \exp(v))) &= u + v.
\end{align*}
Thus $p_1(\log H_0)$ is a subgroup of $V_1$ generated by the
projections of the generators of $H$, and thus in fact $p_1(\log H_0)=p_1(\log
H)$.  Moreover, since $p_1(\log H)$ is a lattice in $V_1$, $H_0$ has
finite covolume.

If $p_j(\log H_0^{(j)})\subset p_j(\log H^{(j)})$, we claim that $p_{j+1}(\log H_0^{({j+1})})\subset p_{j+1}(\log H^{({j+1})})$.  As before, $p_{j+1}(\log H_0^{({j+1})})$ is a subgroup of $V_{j+1}$.
By induction, if $u\in H_0$, $v\in H_0^{(j)}$, there are $u'\in H$, $v'\in H^{(j)}$ such that
$p_1(\log u)=p_1(\log u' )$ and $p_j(\log v )=p_j(\log v')$.
Therefore, again using the Baker-Campbell-Hausdorff formula,
$$p_{j+1}(\log [u,v])=p_{j+1}([\log u,\log v])=\log [p_1(\log u),p_j(\log v)]=p_{j+1}(\log [u',v']),$$
so $p_{j+1}(\log H_0^{({j+1})})\subset p_{j+1}(\log H^{({j+1})}).$
\end{proof}

The techniques in this paper use a lattice in a Lie group together
with a presentation of that lattice, so many theorems require
nilpotent groups that are isomorphic to a lattice in a homogeneous nilpotent
Lie group compatible with the grading.  If $H$ is a finitely-generated
torsion-free nilpotent group of class 2, this is only a mild
restriction; we can show that $H$ has a finite-index subgroup which is
such a lattice.  Consider $H/[H,H]=T\times\Z^d$, where $T$ is torsion.
Then $T[H,H]$ is in the center of $H$.  Let $H_0=(0\times \Z^d)[H,H]$.
Then $H_0$ is finite-index in $H$ and $[H_0,H_0]=[H,H]$, so
$H_0/[H_0,H_0]=\Z^d$.  By Mal\cprime cev's Embedding Theorem, $H_0$ is a
lattice in a connected, 1-connected Lie group $G$ of class 2 with Lie
algebra $\Gamma$.  Any complement of $[\Gamma,\Gamma]$ will give a
grading of $\Gamma$, so since $H_0/[H_0,H_0]=\Z^d$, a choice of $d$
elements of $H_0$ that generate $H_0/[H_0,H_0]$ gives a complement of
$[\Gamma,\Gamma]$ and a grading of $\Gamma$ compatible with the
lattice.

On the other hand, the criterion becomes trickier for groups of
nilpotency class $>2$.  Typically, there are many fewer possible
choices of grading, and it becomes difficult to tell whether a given
presentation gives a lattice compatible with a grading.

\section{Main Theorem}\label{mainthsec}
Gromov\cite{GroAI,GroCC} used the fact that horizontal discs scale well to prove isoperimetric inequalities.
Roughly, if every $(V_1\oplus\dots\oplus V_i)$-horizontal closed curve
$w$ of length $\le 1$ can be filled with a horizontal disc of area at most $c$,
then every horizontal closed curve $w$ of length $\le l$ can be filled
with a disc of area $\sim l^{2i}$ by scaling $w$ to a curve of
length $\le 1$, filling it, then stretching it back again.  In practice,
this criterion is difficult to use, but it suffices to be able to fill
a subset of all horizontal curves of low ``complexity''.  Then, by
expressing the original curve as a combination of ``simple'' curves at
various scales, we can obtain an efficient filling.  For example, to
find fillings for all polygonal curves, it suffices to find
fillings for triangles, then efficiently decompose the original
polygon into triangles.

In \cite{GroCC} and \cite{GroAI}, Gromov used the h-principle for infinitesimally invertible operators to construct horizontal discs.
We will take a combinatorial approach to this technique, using the
fact that a presentation of a group gives fillings of closed curves in
its Cayley graph.  Then, if we can obtain an inequality of the form
$\delta(s_t(r))\lesssim f(t)$ for all the relators $r$ in a
presentation, we get a similar bound $\delta(s_t(w))\lesssim f(t)$
for $w$ any word of bounded length representing the identity.  We then
use these short words as our ``simple'' curves and, by taking
successive approximations, express arbitrary words as combinations of
short words at various scales.  As above, one way to get such a bound
is by filling in relators with a horizontal disc; we will prove
Theorem \ref{mainscale} in the general case, ignoring the source of the
bounds, and then consider some ways of getting such bounds in the
corollaries and examples.

Defining scaling maps on discrete groups and lattices in Lie groups
will be crucial to this approach.  If $v\in V_1$, then $s_t(v)=tv$ and
$s_t(\exp(v))=\exp(v)^t$.  Thus, if $H$ is a lattice compatible with a
grading of $G$, it is generated by $a_i$ such that $a_i=\exp(v_i)$,
where $v_i\in V_1$.  When $t$ is an integer, $s_t(a_i)=a_i^t$
and $s_t$ restricts to a map from $H$ to itself.  This is not,
however, the case in general; if $H$ is a discrete nilpotent group
generated by $a_1,\dots,a_d$, there may not be a homomorphism taking
$a_i$ to $a_i^t$.

We can extend the map from $H$ to $G$ to a map from the Cayley graph of
$H$ into $G$ by mapping in edges as segments of one-parameter
subgroups; then $s_t$ maps the Cayley graph into itself when $t$ is an
integer.  Moreover, if one defines a Carnot-Carath\'{e}odory metric on $G$
as in Section \ref{preliminaries}, $s_t$ acts as a homothety on the
metric, allowing us to compare the metrics $d_{s_t(H)}$ to one another
by comparing them to $d_G$.  We will use this relationship in the
proof of the following theorem.

\begin{thm}\label{mainscale}
Let $G$ be a homogeneous nilpotent Lie group of class $k$ with scaling
automorphism $s_t$.  Let
$$\Gamma=V_1\oplus\dots\oplus V_k$$ 
be the grading on its Lie algebra.  Let $H$ be a lattice in $G$, and
let $H'$ be a nilpotent group with presentation
$\{a_1,\dots,a_d|r_1,\dots,r_s\}$ such that $H$ is a finite quotient $H'/T$ of $H'$ and the presentation
$H=H'/T=\{a_1,\dots,a_d|r_1,\dots,r_s,T\}$ is compatible with the
grading, that is, the images of the $a_i$ lie in $\exp V_1$.

If $\delta_H(s_{2^i}(r_i))\le f(2^i)$ for all $i$, then
$$\delta_{H}(n)\prec \sum_{i=0}^{\lceil\log_2 n\rceil-1} f(2^i)
2^{\lceil\log_2 n\rceil-i}.$$ In particular, if $f(t)\le c
t^\alpha$, then if $\alpha<2$, $\delta(n)$ is subquadratic, so $G$ has a linear isoperimetric inequality\cite{Bow} and is thus $\R$, and if $\alpha\ge 2$, then
$\delta(n)\lesssim n^\alpha$.
\end{thm}
\begin{proof}
By abuse of notation, given a word in the $a_i$, we will consider it
as a word in $H$ or $H'$ interchangeably.

The $r_i$ are assumed to scale efficiently, so given a word $w$
representing the identity in $H$, we would like to construct a disc
filling in $w$ that mostly uses the $r_i$.  If $w$ represents the
identity in $H'$ as well, we can reduce $w$ to the identity using only
the $r_i$.  Then, using the assumption that scalings of the $r_i$ can
be filled efficiently, we find:
\begin{align*}
w&=\prod_{i=1}^{\delta_{H'}(w)} g_i^{-1} r_{p_i}^{\pm1} g_i\\
s_t(w)&=\prod_{i=1}^{\delta_{H'}(w)} s_t(g_i^{-1}) s_t(r_{p_i}^{\pm1}) s_t(g_i)\\
\delta_{H}(s_t(w))&\le\sum_{i=1}^{\delta_{H'}(w)} \delta_{H}(s_t(r_{p_i}^{\pm1}))\\
           &\le \delta_{H'}(w) f(t).
\end{align*}
Thus scalings of words representing the identity in $H'$ can be filled efficiently in $H$.
Note that $s_t(w)$ may not represent the identity in
$H'$, so $\delta_{H'}(s_t(w))$ may not be defined.

Our goal is to break an arbitrary word into pieces which are scalings
of short words representing the identity in $H'$; by the above, these
can be filled efficiently in $H$.  We will proceed by constructing
approximations $w_0,\dots,w_k$ of $w$ and then constructing annuli
between $w_i$ and $w_{i+1}$.  Each $w_i$ will be a word in
$\{a_1^{2^i},\dots,a_d^{2^i}\}$ representing the identity in $H'$.

Let $H_t=s_t(H)=\langle a_1^t,\dots,a_d^t\rangle$ and $H'_t=\langle
a_1^t,\dots,a_d^t\rangle$.  We note that $H'_t/T=H_t$ and that this is
a quasi-isometry; there's a $c_0$ independent of $t$ such that if
$a,b\in H'_t$ project to the same element of $H_t$, then
$d_{H'_t}(a,b)<c_0$.  Let $\rho_t:H\to H_t$ be such that $\rho_t(h)$
is a closest element(in the metric on $G$) of $H_t$ to $h$.  Let
$\rho'_t:H\to H'_t$ be such that $\rho'_t(h)$ is an element of the coset
$\rho_t(h)$.

Let
$$h_{i,j}=\rho'_{2^i}\left(w\left(\left\lfloor\frac{jn}{2^{k-i}}\right\rfloor\right)\right),$$
and let $w_i$ be a closed curve connecting the vertices
$$e=h_{i,0},h_{i,1},\dots,h_{i,2^{k-i}}=e.$$
by shortest paths in $H'_{2^i}$.  In particular, $w_k$ is just a
point and the projection of $w_0$ to $H$ goes through all of the
vertices of $w$.  Note that $w_0\ne w$, since $w$ may not represent the identity in $H'$, but $w_0$ does.  They both represent the identity in $H$, and by applying $cn$ relators of $H$, we can reduce $w_0$ to $w$.

We will need to go between the metrics on all of these groups and their
scalings, which will use the following inequalities:
$$c_1^{-1}d_{G}(h_1,h_2)\le t d_{H_t}(h_1,h_2) \le c_1
d_{G}(h_1,h_2)$$
$$d_{G}(h,\rho_t(h))\le c_2 t$$
which come from the fact that $H$ is quasi-isometric to $G$ and the
scaling automorphism is a homothety of the Carnot-Carath\'{e}odory metric
on $G$, and
$$d_{H_t}(h_1,h_2) \le d_{H'_t}(h_1,h_2)\le d_{H_t}(h_1,h_2)+2c_0.$$
Note in particular that we can go between $d_{H'_t}$ and $d_{H_t}$ with no
loss in the linear constant and between $sd_{H_s}$ and $td_{H_t}$ only
changing the linear constant by $c_1^2$.  

Since $w(\lfloor{jn}/{2^{k-i}}\rfloor)$ and
$w(\lfloor{(j+1)n}/{2^{k-i}}\rfloor)$ are separated in $H$ by at most
$n/2^{k-i}+1$, there is a $c_3$ such that 
$$d_{H'_{2^i}}(h_{i,j},h_{i,j+1})<c_3$$
Thus $w_i$ is a path in $H'_{2^i}$ made up of $2^{k-i}$ segments of length at most $c_3$.

We will fill $w$ by constructing annuli between $w_i$ and $w_{i+1}$.
First, we break this area up into pentagons lying in the Cayley graph
of $s_{2^i}(H')$ with vertices
$$h_{i,2j},h_{i,2j+1},h_{i,2j+2},h_{i+1,j+1},h_{i+1,j}$$ Since the
pentagons lie in the Cayley graph of $s_{2^i}(H')$, they can be filled
efficiently as above.

We bound the length of the perimeter of these pentagons in $s_{2^i}(H')$.  The segments
lying in $w_i$ and $w_{i+1}$ have total length at most $4c_3$.  Since
$$h_{i,2j}=\rho'_{2^i}\left(w\left(\left\lfloor\frac{2jn}{2^{k-i}}\right\rfloor\right)\right)$$
and
$$h_{i+1,j}=\rho'_{2^{i+1}}\left(w\left(\left\lfloor\frac{2jn}{2^{k-i}}\right\rfloor\right)\right),$$
are $\rho'_{2^i}$ and $\rho'_{2^{i+1}}$ of the same point, we find that
there is $c_4$ independent of $i,j$ such that the two remaining segments have length at most
$$d_{H'_{2^i}}(h_{i,2j},h_{i+1,j})\le c_4$$

Thus, each pentagon is a scaling of a short closed curve(of length at most
$4c_3+2c_4$) in $H'$ and thus can be filled in $H$ using at most  
$\delta_{H'}(4c_3+2 c_4)f(2^i)$ relators.  Between $w_i$ and $w_{i+1}$ there are
$2^{k-(i+1)}$ pentagons, and $w_k$ is just a point.  Using $cn$ more
relators to reduce $w$ to $w_0$ gives a filling of $w$ in $H$.
Summing these numbers, we get the required inequality:
$$\delta(n)\le cn+ \sum_{i=0}^{k-1} c_6 2^{k-(i+1)} f(2^i).$$
\end{proof}

Applications of this theorem rely on various ways of finding the bound
$$\delta(s_{2^i}(r_i))\le f(2^i)$$ on the areas of scaled relators.  One
way is to note that commutators like $[x^n,y^n]$ are a product of
$n^2$ conjugates of $[x,y]$:
\begin{cor}\label{quadscale}
With notation as in the previous theorem, if each $r_i$ is
of the form $[x_i,y_i]$, where $x_i$ and $y_i$ are products of commuting
generators, then $G$, $H$, and $H'$ satisfy a quadratic isoperimetric inequality.
\end{cor}
\begin{proof}
We need to show that $\delta_{H}(s_t(r_i))<C t^2$.  If $x$ is a
product of commuting generators, then $s_t(x)=x^t$ as an element of $H'$
and $s_t(x)$ can be reduced to $x^t$ by applying $c t^2$ relators.
Thus, $s_t(r_i)$ can be changed to $[x^t,y^t]$ by applying $4c t^2$
relators and $[x^t,y^t]$ can be reduced to the identity by commuting
$x$ and $y$ $t^2$ times.
\end{proof}

Similarly, commutators of the form $[x_1^n,\dots,x_k^n]$ are products
of $n^k$ conjugates of $[x_1,\dots,x_k]$, and a similar theorem
applies to groups which can be presented with relations of the form
$[x_1,\dots,x_k]$.  In particular, this gives a bound
$\delta(n)\lesssim n^{c+1}$ for free nilpotent groups of class $c$.

Though the corollary applies in the case that $H'$ is a nilpotent group of
class $>2$, I have been unable to find an example of such a group for
which the hypotheses hold, though the class 3 group considered in Section
\ref{appssec} comes close.  In the next section, we will consider
some examples of groups of class 2 to which the corollary applies.

In the above corollary, the same bounds could have been obtained using
differentiable methods; if $x$ and $y$ are as in the theorem, it is easy to
obtain a horizontal disc filling $[x,y]$ and thus a quadratic bound on
the area of a scaling.  The following corollary gives an example of a
bound where the combinatorial aspects are more essential.

\begin{cor}\label{quotdehn}
Let $G$ be a homogeneous nilpotent Lie group of class $k$ with scaling
automorphism $s_t$.  Let
$$\Gamma=V_1\oplus\dots\oplus V_k$$ 
be the grading on its Lie algebra.  Let $H$ be a lattice in $G$ compatible with the grading.

Assume that $H$ satisfies the
isoperimetric inequality $\delta_{H}(n)\prec n^\alpha$.  If
$r_1,\dots,r_s\in H$ are such that $\log r_i\in V_{j_i}$, then the
Dehn function $\delta_{H/\{r_1,\dots,r_s\}}$ is bounded by
$$
\delta_{H/\{r_1,\dots,r_s\}}(n)\prec \begin{cases}n^\alpha & \text{if $\alpha>j$}\\
               n^\alpha \log n & \text{if $\alpha=j$}\\
               n^j & \text{if $\alpha<j$}
				\end{cases},
$$
where $j=\max\{j_i\}$.
\end{cor}
\begin{proof}
We will apply Theorem \ref{mainscale} to $\hat{H'}= H/\{r_1,\dots,r_s\}$.  Let
$\mathfrak{r}\subset\Gamma$ be the Lie ideal generated by the $\log r_i$.
Since $\mathfrak{r}$ is a graded ideal, $\Gamma/\mathfrak{r}$ determines
a homogeneous simply connected Lie group.  We have a map from $H/\{r_1,\dots,r_s\}\to
\Gamma/\mathfrak{r}$; we claim that the kernel of this map is the set of
torsion elements of $\hat{H'}$ and the image is a lattice in
$\Gamma/\mathfrak{r}$.  Clearly, the kernel contains the torsion elements;
moreover $\{r_1,\dots,r_s\}$ generate a lattice in $\mathfrak{r}$ and thus
generate a subgroup of finite index in $H \cup
\mathfrak{r}$. Therefore, the kernel is finite, and so all of its elements
are torsion.  Similarly, since $\{r_1,\dots,r_s\}$ generate a lattice
in $\mathfrak{r}$, the image is a lattice in $\Gamma/\mathfrak{r}$.  Call this
image $\hat{H}$.  $\hat{H}$ and $\hat{H'}$ satisfy the hypotheses of
Theorem \ref{mainscale}; we need to find a bound on the area of scaled
relators.  Given a presentation for $H$, we can get a presentation
for $\hat{H'}$ by appending the $r_i$.  It thus suffices to find
bounds on the area of scalings of words representing the $r_i$.

Consider an element $r\in H$ with $\log r \in V_j$.  Let $w$ be a word in
$H$ representing $r$.  We want to bound
$\delta_{\hat{H}}(s_{2^n}(w))$.  But since $\log r\in V_j$,
$s_t(r)=r^{t^j}$ as elements of $H$ so $s_{2^n}(w)s_{2^{n-1}}(w)^{-2^j}$ represents the
identity in $H$.  We can thus reduce $s_{2^n}(w)$ to $2^j$ copies of
$s_{2^{n-1}}(w)$ by applying $$\delta_{H}\left(2^n l(w)+2^j\cdot
2^{n-1} l(w)\right)=\delta_{H}\left(2^{n}(1+2^{j-1}) l(w)\right)$$ relators.

Repeating this process $n$ times, we can reduce $s_{2^n}(w)$ to
$2^{nj}$ copies of $w$ by using
\begin{align*}
\sum_{i=0}^n 2^{(n-i)j}\delta_{H}\left(2^{i}(1+2^{j-1}) l(w)\right)&<2^{nj}c_0\sum_{i=0}^n 2^{-ij}2^{i\alpha}\\
&<\begin{cases}c_1 2^{n\alpha} & \text{if $\alpha>j$}\\
               c_1 n 2^{n\alpha} & \text{if $\alpha=j$}\\
               c_1 2^{n j} & \text{if $\alpha<j$}
\end{cases}
\end{align*}
relators, where the $c_i$ are constants independent of $n$.  This
final word takes another $2^{nj}$ relators to reduce, giving the required bound.
\end{proof}

This gives a proof of a $n^{c+1}$ bound for homogeneous nilpotent Lie
groups with rational structure constants; such a group is a quotient
of a free nilpotent group by a graded ideal $\mathfrak{r}$.  By the
remark after Corollary \ref{quadscale}, the free nilpotent group
satisfies the isoperimetric inequality $\delta(n)\lesssim n^{c+1}$.
Corollary \ref{quotdehn} applies with $r_i$ generating $\mathfrak{r}$;
we can choose such generators so that the corollary gives an
isoperimetric inequality $\delta(n)\lesssim n^{c+1}$.  This result was first proved by Pittet for all homogeneous nilpotent Lie groups in \cite{PitHom}.

\section{Applications and Examples} \label{appssec}
The hypotheses for Corollary \ref{quadscale} are very restrictive and
may exclude many nilpotent groups with quadratic Dehn functions.
Nonetheless, many interesting nilpotent groups can be given such a
presentation, and one source of such examples is central powers.  Many
methods are available for working with central powers of nilpotent
groups; Mark Sapir informs me that he and Ol'shanskii used the
methods of \cite{SapOls} to prove some of the following results.
Given $G$ a nilpotent group of class $k$, define
$G^{\times_{G^{(k)}}n}$ to be $G^n$ with the copies of $G^{(k)}$
identified.  Then, for example, the $(2n+1)$ dimensional Heisenberg
group is the $n$th central power of $H_3$, the 3-dimensional
Heisenberg group.  When $n\ge 2$, these groups have been shown to have
a quadratic isoperimetric inequality; Corollary
\ref{quadscale} gives an alternate proof of this fact as follows:

$H_3$ can be given the presentation
$$H_3=\{a_1,a_2,c|[a_1,a_2]=a_3,[c,a_i]=1\}.$$ 
$H_5$ can be given the presentation
$$H_5=\{a_1,a_2,b_1,b_2,c|[a_1,a_2]=[b_1,b_2]=c,[a_i,b_i]=[a_i,c]=[b_i,c]=1\}.$$
However, many of these relators are redundant; the relations
$[a_i,c]=1$ follow from the fact that $[a_i,c]=[a_i,[b_1,b_2]]=1$, and we can eliminate $c$ to get the presentation
$$H_5=\{a_1,a_2,b_1,b_2|[a_1,a_2][b_1,b_2]^{-1}=1,[a_i,b_i]=1\}.$$
Furthermore, since the $a_i$ and $b_i$ commute, we can interleave $[a_1,a_2]$ and $[b_1,b_2]^{-1}=[b_2,b_1]$ to get:
\begin{align*}
[a_1,a_2][b_1,b_2]^{-1}&=a_1a_2a_1^{-1}a_2^{-1}b_2b_1b_2^{-1}b_1^{-1}\\
                       &=a_1b_2a_2b_1b_2^{-1}a_1^{-1}b_1^{-1}a_2^{-1}\\
                       &=[a_1b_2,a_2b_1],
\end{align*}
so
$$H_5=\{a_1,a_2,b_1,b_2|[a_1b_2,a_2b_1]=1,[a_i,b_i]=1\}.$$
We can then apply Corollary \ref{quadscale} to show that $H_5$ has a
quadratic isoperimetric inequality.

This technique can be generalized to:
\begin{cor}\label{centcor}
Let $H$ be a nilpotent group of class 2 with a presentation of the
form $$H=\{a_1,\dots,a_p|r_1,\dots,r_s,[a_i,[a_j,a_k]]=1\}.$$
Let $Z=[H,H]$ and assume that the $r_i$ are
of the form $\prod_{j=1}^{l_i}[x_j,y_j]=1$, where the $x_j$ and $y_j$
are products of commuting $a_i$.  Let $n\ge \max\{2,l_i\}$  Then the central
product of $n$ copies of $H$(with the $i$th copy generated by $a_{i1},\dots,a_{ip}$):
\begin{align}
\nonumber H^{\times_Z n}=&\{a_{11},\dots,a_{np}|\\
\label{origrels}   &r_{i1},\dots,r_{is}\text{\qquad for $1\le i\le n$}\\
\label{nilprels}   &[a_{ij},[a_{ik},a_{il}]]=1\text{\qquad for all $i,j,k,l$}\\
\label{productrels}&[a_{ij},a_{kl}]=1\text{\qquad for $i\ne k$}\\
\label{centerrels} &[a_{ik},a_{il}]=[a_{jk},a_{jl}]\text{\qquad for all $i,j,k,l$}\}
\end{align}
satisfies a quadratic isoperimetric inequality.
\end{cor}
\begin{proof}
Note that the presentation given is indeed the central product as
defined before.  The relators in (\ref{origrels}) and (\ref{nilprels}) are
$n$ copies of the presentation of $H$, the relators in
(\ref{productrels}) imply that the group is a quotient of the direct
product of those $n$ copies, and the relators in (\ref{centerrels})
identify the copies of $Z$.

Moreover, $H/[H,H]=\Z^{p}$ and the projections of the $a_{i}$
generate $H/[H,H]$.  Thus, by the discussion in Section
\ref{preliminaries}, if $T$ is the subgroup of torsion elements of $H$, $H/T$ is a lattice in a homogeneous Lie group compatible
with the grading.

We want to express $H^{\times_Z n}$ in the form used by Corollary
\ref{quadscale}; to do that, we need to modify the relators in
(\ref{origrels}) and (\ref{centerrels}) and eliminate those in
(\ref{nilprels}).

The relators in (\ref{nilprels}) follow from the relations in
(\ref{productrels}) and
(\ref{centerrels}).  By the relations in (\ref{centerrels}),
$$[a_{ik},[a_{il},a_{im}]]=[a_{ik},[a_{jl},a_{jm}]],$$ which is
the identity, since the relations in (\ref{productrels}) imply that
$a_{ik}$ commutes with any word in the $a_{j,\cdot}$.

Using the relators in (\ref{productrels}) as for the Heisenberg group,
we can change the relators in (\ref{centerrels}) into commutators:
$$[a_{ik},a_{il}][a_{jk},a_{jl}]^{-1}=[a_{ik}a_{jl},a_{il}a_{jk}]$$

A similar process makes the relators in (\ref{origrels}) into
commutators.  Let $r_{ij}$ be the relator
$\prod_{k=1}^{l_j}[x_{ij},y_{ij}]=1$, where $x_{ij}$ and $y_{ij}$
are the words $x_j$ and $y_j$ written in the $a_{i,\cdot}$.  By using
the relations in (\ref{nilprels}), we can express $[x_{ij},y_{ij}]$ as a product of
commutators of the $a_{i,\cdot}$.  The relations in (\ref{centerrels}) let us replace these commutators with commutators of the $a_{i',\cdot}$ for any $i'$, so we can replace $[x_{ij},y_{ij}]$
with $[x_{i',j},y_{i',j}]$.  Replace
$\prod_{k=1}^{l_j}[x_{ij},y_{ij}]$ with
$\prod_{k=1}^{l_j}[x_{kj},y_{kj}]$.  All the commutators are now
written in different sets of generators, so as before, we can
interleave them to obtain a relation
$$\left[\prod_{k=1}^{l_j}x_{kj},\prod_{k=1}^{l_j}y_{kj}\right]=1$$
which is of the form required for Corollary \ref{quadscale}.
\end{proof}

Presentations of this form are common.  If $H$ is a finitely generated
nilpotent group of class 2 such that $H/[H,H]=\Z^d$, then it has a
presentation of the form required by Corollary \ref{centcor}, and any
torsion-free nilpotent group has a finite index subgroup with
$H/[H,H]=\Z^d$.  Thus any finitely generated nilpotent
group of class 2 has a central power with a quadratic Dehn function.

The power required depends on the group, but many nilpotent groups
only require a low power.  If $r$ is a relator of the form
$\prod_{j=1}^{l}[x_j,y_j]=1$, where the $x_j$ and $y_j$ are products
of commuting $a_i$, we will call $l$ the complexity of $r$.  Then if
$H$ has a presentation with relations $[a_i,[a_j,a_k]]$ and relations
of  complexity at most $n$, then the $\max\{n,2\}$th central power of $H$
has quadratic Dehn function.  Many groups have presentations using
simple relators.  One example is free nilpotent groups; if $G$ is the
free nilpotent group of class $2$ with $d$ generators, then it has a
presentation $$G=\{a_1,\dots,a_d|[a_i,[a_j,a_k]]=1\}$$ and $G\times_Z
G$ has quadratic Dehn function.

\label{octonexamplesec}
Certain generalizations of the Heisenberg group can be presented using
only relations of length 2.  If $k$ is the complex numbers,
quaternions, or octonions, considered as a real vector space, we can
define a nilpotent Lie algebra $\Gamma$ of class 2 where
$\Gamma/[\Gamma,\Gamma]=k$ and $[\Gamma,\Gamma]$ is the subspace of imaginary
elements of $k$ by letting $[v,w]=\imag(v\bar{w})\in [\Gamma,\Gamma]$.
The corresponding Lie groups have lattices whose elements correspond to elements
of $k$ with integer coefficients.  These groups show up as cusp
subgroups in groups acting on negatively curved symmetric spaces.  In
the case that $k=\C$, this is the three-dimensional Heisenberg group
and its central powers are the higher-dimensional Heisenberg groups.

These groups can be presented with relations of length 2; in the
case of the quaternions, for instance, a typical relator looks like:
$$[i,j]=[1,k],$$ 
where $1,i,j$ and $k$ represent generators of $G$ corresponding to
unit quaternions.  Thus, central powers of these groups have lattices
with quadratic Dehn functions and thus the corresponding Lie groups
have quadratic Dehn functions as well.

In general, we have the following bound:
\begin{cor}
If $G$ is a nilpotent Lie group of class 2 with a basis with rational
structure constants, then $G^{\times_Z n}$ has a quadratic
isoperimetric inequality when $n\ge\max\{2,\dim [G,G]\}$.  In
particular, if $H$ is a finitely-generated nilpotent group of class 2,
then $H^{\times_Z n}$ has a quadratic isoperimetric inequality when
$n\ge\max\{2,\rk [H,H]\}$.
\end{cor}
\begin{proof}
If $G$ is such a group, it contains lattices which are compatible with
some grading.  We need to show that some lattice has, up to finite
index, a presentation with relators of complexity $n$.  Let
$a_1,\dots,a_d$ be part of a basis of $G$ that has rational structure
constants and generates $G/[G,G]$.  Then $[G,G]=\R^{\binom{d}{2}} /
R$ for some $R\subset \R^{\binom{d}{2}}$, where $\R^{\binom{d}{2}}$ has a
basis consisting of elements $[a_i,a_j]$.  We can find a basis of $R$
consisting of vectors which are each in the span of at most $n$ basis
vectors and have integer coordinates in the basis.  These correspond
to relators of complexity at most $n$ in the presentation of a group
$H$ with a finite-index quotient $H/T$ which is a lattice in $G$.  
\end{proof}

Smaller central powers still have strong isoperimetric
inequalities:
\begin{cor}\label{centcorlog}
If $H$ is a nilpotent group of class 2, and $Z=[H,H]$,
then $H^{\times_Z p}$ satisfies the isoperimetric inequality
$\delta(l)\lesssim l^2 \log l$. 
\end{cor}
\begin{proof} 
Up to finite index, $H/[H,H]=\Z^d$.  Then $H$ is a quotient of the free
nilpotent group $N=F_d/F_d^{(3)}$ of class 2 on $d$ generators by some
relators contained in $[N,N]$.  Similarly, $H^{\times_Z
  p}$ is a quotient of $N^{\times_{[N,N]} p}$ by relators contained in
$[N^{\times_{[N,N]} p},N^{\times_{[N,N]} p}]$.  We can thus apply
Corollary \ref{quotdehn}. 

By the remark after Corollary \ref{centcor}, $[N^{\times_{[N,N]}
p},N^{\times_{[N,N]} p}]$ satisfies a quadratic isoperimetric
inequality, so $H^{\times_Z p}$ satisfies the isoperimetric inequality
$\delta(l)\lesssim l^2 \log l$. 
\end{proof} 

It is an interesting question whether these groups necessarily have
quadratic Dehn functions.  All of the nilpotent groups for which we
know quadratic isoperimetric inequalities have such inequalities
because curves can be filled by horizontal
discs\cite{GroCC,Allcock}; in Section \ref{vscentiso} we will
give an example of a group where Corollary \ref{centcorlog} gives an $n^2 \ln n$ upper bound, but
there are curves which cannot be filled by such discs.

Finally, Theorem \ref{mainscale} allows us to give an example of a nilpotent Lie
group of class 3 and rank 8 satisfying a quadratic isoperimetric
inequality.  This group has Lie algebra $\Gamma$ with the grading:
$$\Gamma=\langle a,b,c,d,e \rangle\oplus\langle f,g
\rangle\oplus\langle h\rangle.$$
and bracket given by
\begin{align*}
[a,b]&=[c,d]=f\\
[b,c]&=[d,e]=g\\
[b,f]&=[g,d]=h
\end{align*}
with all other brackets 0.

This is generic among homogeneous nilpotent Lie algebras
with a grading having these dimensions.  A class 2 nilpotent Lie
algebra $\Phi$ can be described by the skew-symmetric bilinear form
$[\cdot,\cdot]:\Phi/[\Phi,\Phi]\to[\Phi,\Phi]$; a generic pair of
skew-symmetric forms on $\R^5$ give a Lie algebra isomorphic to $\Gamma/\Gamma^{(3)}$(see \cite{Tho} or \cite{LancRod} for the classification of pairs of skew-symmetric forms).
By calculating cohomology, one can show that $\Gamma/\Gamma^{(3)}$ has a unique central extension to a
homogeneous class 3 nilpotent Lie algebra.

The proof that $G$ has a quadratic Dehn function involves a lengthy
combinatorial group theory calculation, given in Appendix
\ref{combgp}, but we sketch the proof here:
$\exp(a),\exp(b),\exp(c),\exp(d),\exp(e)$ generate a lattice $H$.  By abuse
of notation, we will consider $a,b,c,d,e$ to be elements in this
lattice.  The relations
$$[a,c],[a,d],[a,e],[b,d],[b,e],[c,e],[ac^{-1}e,bd],[ace,bd^{-1}]$$
hold in $H$.  Since they are of the form required by Theorem
\ref{mainscale}, they can be filled by horizontal discs.  These
relations do not present $H$, but we show that any word representing
the identity in $H$ can be filled by a horizontal disc composed of
scalings and distortions of those discs.  Thus, given a presentation,
we can fill its relators with horizontal discs and Theorem
\ref{mainscale} will give us the required quadratic isoperimetric
inequality.

\section{Conclusion}
Most of the isoperimetric inequalities proved here rely on filling
curves with horizontal discs.  In fact, all
isoperimetric inequalities for lattices in nilpotent groups arise from
fillings by almost-horizontal discs:
\begin{thm}\label{almosthoriz}
Let $G$ be a homogeneous nilpotent Lie group of class $k$ with scaling
automorphism $s_t$.  Let
$$\Gamma=V_1\oplus\dots\oplus V_k$$
be the grading on its Lie algebra.  Let $H\subset G$ be a
lattice in $G$ compatible with the grading.  

Fix a metric on $P\Gamma$ the projective tangent space of $G$ and call
a piecewise $C^2$ map whose tangent vectors are non-zero in the
interior of a piece and $\epsilon$-close to $V_1$
{\em$\epsilon$-horizontal}.

If $\delta(l)$ is the Dehn function of $H$ and $\gamma$ is a word in
$H$, there are $c_i$ such that for $t$ sufficiently large, $\gamma$ can be filled with
a $c_0/t$-horizontal disc of area at most $c_1 \delta(c_2 t)/t^2$.
\end{thm}
\begin{proof}
If $V_1$ is 1-dimensional, $V=\R$ and there is nothing to prove.
Otherwise, $\epsilon$-horizontality is an open $\Diff$-invariant
differential relation for all $\epsilon>0$ and thus a microflexible
differential relation.  In particular, by the arguments in
\cite{GroCC}, any horizontal curve can be filled by an $\epsilon$-horizontal disc.

For $\epsilon$ sufficiently small, we can choose $c_0$ such that if
$f$ is an $\epsilon$-horizontal disc, then for all $t$, $s_{1/t}\circ
f$ is $c_0/t$-horizontal.

Fill each relator of $H$ with an $\epsilon$-horizontal disc and let
$c_1$ be the maximum area of such a disc.  Then when $t$ is an
integer, $s_t(\gamma)$ is a word that can be filled by an
$\epsilon$-horizontal disc of area at most $c_1 \delta(l(\gamma)t)$.
Applying $s_{1/t}$ to this disc gives us a $c_0/t$-horizontal filling
of $\gamma$ with area at most $c_1' \delta(l(\gamma)t)/t^2.$
\end{proof}

It is thus interesting to consider groups where some curves cannot be
filled by horizontal discs, as this failure may lead to a lower bound
for the Dehn function.  One source of such examples comes from taking
quotients of groups as in Corollary \ref{quotdehn}; it's possible to
take quotients by relators that satisfy the conditions of the theorem
but which cannot be filled by horizontal discs.

\label{vscentiso}
One such example, due to Mark Sapir, can be
constructed as follows.  Let $\Phi$ be the free nilpotent Lie algebra
of class 2 with $\Phi/[\Phi,\Phi]\cong \R^{10}$ and take
$Z=[\Phi,\Phi]$.  Given a Lie algebra, we can consider it as a simply
connected Lie group with multiplication given by the
Baker-Campbell-Hausdorff formula.  Then $\Phi\times_Z \Phi$ has a
quadratic Dehn function.  Take the quotient $$\Phi\times_Z
\Phi/[a_1,a_2][a_3,a_4]\dots[a_9,a_{10}],$$ where the $a_i$ are generators
of $\Phi/[\Phi,\Phi]$.  By Corollary \ref{quotdehn}, this satisfies
the isoperimetric inequality $\delta(l)\lesssim l^2\log l$.  Let
$\Psi$ be its Lie algebra.  Then $\Phi$ and $\Psi$ have gradings
$\Phi=V_1\oplus V_2$ and $\Psi=W_1\oplus W_2$, with $V_1\cong
W_1=\langle a_1,\dots,a_{10}\rangle$.  Any $V_1$-horizontal disc in
$\Phi$ projects to a disc in $V_1$, and if $v,w$ span the tangent
plane at a point, then the fact that the disc is $C^2$ implies that
$[v,w]=0$.  Conversely, any disc whose tangent spaces satisfy this
condition can be lifted to a horizontal disc in $\Phi$, unique up to
translation by an element of $V_2$.  However, if $v,w\in V_1$, then
$[v,w]_\Phi=0$ iff $[v,w]_\Psi=0$, so the set of such projections of
horizontal discs is the same for $\Phi$ and $\Psi$.  Thus, since
$[a_1,a_2]\dots[a_9,a_{10}]$ doesn't represent the identity in
$\Phi\times_Z \Phi$, there is no horizontal disc in $\Phi\times_Z
\Phi$ filling it, and thus no horizontal disc in $\Phi\times_Z
\Phi/[a_1,a_2][a_3,a_4]\dots[a_9,a_{10}]$ filling it.

Another example is the nilpotent group based on the octonions
described in Section \ref{octonexamplesec}; on the one hand, the
centralized isoperimetric function is quadratic(Pittet's original calculation in \cite{Pittet} had a sign error, corrected in \cite{LeuzPit}), so
the best known lower bound for the isoperimetric function is
quadratic, but there are no horizontal discs.  As above, the
generators of any tangent plane to a horizontal disc must commute, but
$[v,w]=\imag(v\bar{w})$, which is $0$ if and only if $v$ and $w$ are
multiples of one another or $0$.

A lack of horizontal discs, however, does not imply a lack of
almost-horizontal discs or a lower bound for the isoperimetric
inequality.  The h-principle for open differential relations implies
that almost-horizontal discs always exist and does not give the lower
bound on their size required by Theorem \ref{almosthoriz}.  Finding a
quantitative version of this h-principle which gives such bounds would
be a powerful tool for dealing with the Dehn functions of
nilpotent groups.

\appendix
\section{A nilpotent group of class 3 with quadratic isoperimetric inequality}\label{combgp}
In this appendix, we prove that the class 3 nilpotent group described
in Section \ref{appssec} satisfies a quadratic isoperimetric
inequality.

We recall the definition:

Consider the Lie algebra $\Gamma$ with the grading:
$$\Gamma=\langle a,b,c,d,e \rangle\oplus\langle f,g
\rangle\oplus\langle h\rangle.$$
and bracket given by
\begin{align*}
[a,b]&=[c,d]=f\\
[b,c]&=[d,e]=g\\
[b,f]&=[g,d]=h
\end{align*}
with all other brackets 0.

We consider the simply-connected Lie group $G$ with Lie algebra
$\Gamma$.  Then since $a$, $b$, $c$, $d$, $e$, $f$, $g$, $h$ is a basis with rational
structure constants, $\exp(a)$, $\exp(b)$, $\exp(c)$, $\exp(d)$, $\exp(e)$
generate a lattice $H$ in $G$; by an abuse of notation, we will consider
$a,b,c,d,e$ to be generators of $G$.  

In $\Gamma$, we have the identities:
$$[a,c]=[a,d]=[a,e]=[b,d]=[b,e]=[c,e]=[a+c+e,b-d]=[a-c+e,b+d]=0.$$
Thus in $G$, the words
$$[a,c],[a,d],[a,e],[b,d],[b,e],[c,e],[ace,bd^{-1}],[ac^{-1}e,bd]$$
can be filled with horizontal discs.  Using these words as relators and inserting and deleting them from words corresponds to
gluing together horizontal discs.  Thus we aim to show that any word in the
alphabet $a,b,c,d,e$ that represents the identity in $G$ can be
reduced to the identity by using these words and others like them.

These others will come from automorphisms of $\Gamma$.
We first note that $\Gamma$ has many automorphisms; given real numbers
$t_a,t_b,t_c,t_d,t_e\ne 0$ such that
\begin{align*}
t_at_b&=t_ct_d\\
t_bt_c&=t_dt_e.
\end{align*}
then 
\begin{align*}
a&\mapsto t_a a &b&\mapsto t_b b\\
c&\mapsto t_c c &d&\mapsto t_d d\\
e&\mapsto t_e e &f&\mapsto t_c t_d f\\
g&\mapsto t_b t_c g &h&\mapsto t_b t_c t_d h
\end{align*}
is an automorphism of $\Gamma$.  In particular, any nonzero choice of three of the coefficients gives an automorphism.  Similarly,
\begin{align*}
a&\mapsto e &b&\mapsto d\\
c&\mapsto c &d&\mapsto b\\
e&\mapsto a &f&\mapsto -g\\
g&\mapsto -f &h&\mapsto h
\end{align*}
is an automorphism.  Both of these preserve $V_1=\langle
a,b,c,d,e\rangle$, so the images of horizontal discs remain
horizontal.  

We will essentially be doing a combinatorial group theory calculation,
showing that
$$G=\{a,b,c,d,e|[a,c],[a,d],[a,e],[b,d],[b,e],[c,e],[ace,bd^{-1}],[ac^{-1}e,bd],\dots\},$$
where the dots represent all images of the relators under these
automorphisms.  Any word that can be reduced to the identity using
these relators corresponds to a curve in $\Gamma$ that can be filled
with a horizontal disc.  Note, however, that in general, the image of
a word under an automorphism is a curve in the Lie group rather than
another word.  We will ignore this, referring to, for example, a word
in $a$ and $b$ rather than a curve whose tangent vectors lie in the
$a$ and $b$ directions.

In the following calculations, we will denote the inverses of
$a,$ $b,$ $c,$ $d,$ $e$ by $A,$ $B,$ $C,$ $D,$ $E$.  In addition, $[x,y]$ will denote
$xyx^{-1}y^{-1}$ and $0$ will denote the trivial word.  By abuse of
notation, $w=w'$ will denote that there is a horizontal annulus
between $w$ and $w'$, i.e., that $w$ can be reduced to $w'$ using our
relators or that $w{w'}^{-1}$ can be reduced to the identity using the
relators.  Words may be written in a fixed-pitch font for ease of
reading; this has no mathematical significance.

It is convenient to state a version of the Jacobi identity:
$$[x,[y,z]]=g_1^{-1}[y^{-1},[z^{-1},x^{-1}]]g_1g_2^{-1}[[x,y^{-1}],z^{-1}]g_2$$
In particular, if two of the commutators vanish, so does the third, and modulo $H^{(4)}$, this is just
$$[x,[y,z]]=[y,[x,z]][z,[y,x]].$$

\begin{lemma}
$[c,[b,a]]$ can be filled with a horizontal disc.
\end{lemma}
\begin{proof}
We start with:
$$\mathtt{0=[Bd,ace]cc[aCe,BD]CC=BdacebDACEccaCeBDAcEbdCC}.$$
Moving pairs of $a$'s and $b$'s together cancels them:
\begin{align*}
&\mathtt{BdacebDBDAcEbdCC}\\
&\mathtt{BdaceDDAcEbdCC}\\
&\mathtt{BdceDDcEbdCC}
\end{align*}
Conjugating by $b$, we can collect the occurrences of $b$ into a commutator:
\begin{align*}
&\mathtt{dceDDcEdbCCB}\\
&\mathtt{dceDDcEdCC[cc,b]}
\end{align*}

Thus we can express $[cc,b]$ as a word in $c,d,$ and $e$.  $c,d$, and
$e$ commute with $a$, so $[a,[cc,b]]$ can be filled with a horizontal
disc.  Applying an automorphism to $[a,[cc,b]]$ gives $[a,[c,b]]$.  In
fact, by applying automorphisms, these are zero when $a$, $b$, and $c$
are replaced by powers of themselves.  Since $[a,[c,b]]$ and
$[b,[c,a]]$ can be reduced to 0, we can reduce $[c,[b,a]]$ as well.
\end{proof}

Similarly, by applying various automorphisms, we can show that
$[e^{n_e},[d^{n_d},c^{n_c}]]$ and any of its permutations can be
filled by such a disc for any choice of $n$'s.  

Now consider
$$\mathtt{0=[bd,AcE][ace,bD]=bdAcEBDaCeacebDACEBd.}$$
Several letters cancel to produce
$$\mathtt{dAcEBDaeaebDACEd}$$
We conjugate and collect $a$'s and $b$'s to get
\begin{align*}
&\mathtt{BaabAAEDeeDCEddc}\\
&\mathtt{[B,aa]EDeeDECddc}
\end{align*}
In particular, as in the lemma, $[a,[a,b]]\sim0$.  Moreover, applying
an automorphism, we find that $[a^{n_a},b^{n_b}]$ can be expressed as
a word in $c,d$, and $e$, so it commutes with any power of $a$.

Applying an automorphism to this that takes each generator to its inverse, we get:
$$\mathtt{[b,AA]edEEdecDDC}$$
and combining these, we get
\begin{align*}
&\mathtt{cDDC[b,AA]edEEdeEDeeDECddc[B,aa]}\\
&\mathtt{cDDC[b,AA]Cddc[B,aa]}
\end{align*}
$c$ and $d$ both commute with $[b,a]$, by the Lemma and because $d$
commutes with $a$ and $b$, so we find that 
\begin{align*}
&\mathtt{cDDC[b,AA]Cddc[B,aa]}\\
&\mathtt{[dd,cc][b,AA][B,aa]}
\end{align*}

Since $[a,[a,b]]=0$, this can be reduced to 
$$\mathtt{[dd,cc][b,AA][AA,B]=0}$$ or $$\mathtt{[dd,cc]b[AA,BB]B=0}$$
or, applying an automorphism taking $a\to A$ and $b\to B$,
$$\mathtt{[dd,cc]B[aa,bb]b=0}.$$

Thus $[dd,cc]$ can be reduced to a word in $a$ and $b$, and by
applying automorphisms, $[d^{n_d},c^{n_c}]$ can be reduced to a word
in $a$ and $b$.  So $[d^{n_d},[d^{n_d'},c^{n_c}]]\sim 0$, since $d$
commutes with $a$ and $b$, and similarly for
$[b^{n_b},[b^{n_b'},c^{n_c}]]$.  Our goal is to show that all
4-element commutators can reduced to the identity.  So far,
we've shown that:
\begin{align*}
&\mathtt{[a,[a,b]]=[c,[a,b]]=[d,[a,b]]=[e,[a,b]]=0}\\
&\mathtt{[a,[b,c]]=[b,[b,c]]=[e,[b,c]]=0}\\
&\mathtt{[a,[c,d]]=[d,[c,d]]=[e,[c,d]]=0}\\
&\mathtt{[a,[d,e]]=[b,[d,e]]=[c,[d,e]]=[e,[d,e]]=0}.
\end{align*}
In fact, we have shown that all of these relations hold when $a,b,c,d,e$ are replaced by powers of themselves.

The remaining 3-element commutators are $[c,[c,d]$, $[c,[b,c]]$, $[b,[a,b]]$, $[d,[d,e]]$, $[b,[c,d]]$, and $[d,[b,c]]$.  
\begin{itemize}
\item
$[c,[c,d]]$ and $[c,[b,c]]$:

It suffices to show that $[cc,[cc,dd]]= 0$.  
\begin{align*}
&\mathtt{[cc,[cc,dd]]=[cc,B[aa,bb]b]}\\
&\mathtt{[cc,[cc,dd]]=B[bccB,[aa,bb]]b}\\
&\mathtt{[cc,[cc,dd]]=B[[b,cc]cc,[aa,bb]]b}
\end{align*}

Since $c$ commutes with $[a,b]$ this can be reduced to
$$\mathtt{[cc,[cc,dd]]=B[[b,cc],[aa,bb]]b}$$
but $[b,c]$ commutes with $a$ and $b$, so 
$$\mathtt{[cc,[cc,dd]]=0}$$

Similarly, $[c,[b,c]]=0$.  
\item
$[b,[a,b]]$ and $[d,[d,e]]$:

These can be reduced to conjugates of $[b,[c,d]]$ and $[d,[b,c]]$:
\begin{align}
\nonumber&\mathtt{[bb,[cc,dd]]=[bb,B[aa,bb]b]}\\
\nonumber&\mathtt{[bb,[cc,dd]]=B[bb,[aa,bb]]b}\\
\label{bcdbab}&\mathtt{[bb,[cc,dd]]=B[bb,[aa,bb]]b}
\end{align}
and similarly, $[dd,[bb,cc]]=D[dd,[dd,ee]]d$.
\item
$[b,[c,d]]$ and $[d,[b,c]]$:

We want to show that these are central; this will show that all 4-element commutators can be filled with horizontal discs and lead to our conclusion.  From above, we have that $[dd,[bb,cc]]=D[dd,[dd,ee]]d$, so $[d,[b,c]]$ commutes with $a$ and $b$.  

So consider:
\begin{align*}
\mathtt{[dd,[bb,cc]]}&\mathtt{=ddbbccBBCCDDccbbCCBB}\\
\mathtt{[dd,[bb,cc]]}&\mathtt{=bbcc[[CC,dd],BB]CCBB}
\end{align*}
Since $b$ commutes with $[d,[b,c]]$, 
$$\mathtt{[dd,[bb,cc]]=cc[[CC,dd],BB]CC}$$
Using (\ref{bcdbab}), we find that
\begin{align*}
\mathtt{[dd,[bb,cc]]}&\mathtt{=ccb[[aa,BB],BB]BCC}\\
\mathtt{[dd,[bb,cc]]}&\mathtt{=bcc[CC,B] [[aa,BB],BB] [B,CC]CCB}
\end{align*}
$[c,b]$ commutes with $a$ and $b$, so 
\begin{align*}
\mathtt{B[dd,[bb,cc]]b}&\mathtt{=cc[[aa,BB],BB]CC}\\
\mathtt{[dd,[bb,cc]]}&\mathtt{=cc[[aa,BB],BB]CC}
\end{align*}
$c$ commutes with $[a,b]$, so
\begin{align*}
\mathtt{[dd,[bb,cc]]}&\mathtt{=[[aa,BB],ccBBCC]}\\
\mathtt{[dd,[bb,cc]]}&\mathtt{=[[aa,BB],BB[cc,BB]]}
\end{align*}
Again, $[c,b]$ commutes with $a$ and $b$, so
$$\mathtt{[dd,[bb,cc]]=[[aa,BB],BB]}$$

Thus $[d,[b,c]]$ can be reduced to a word in $a$ and $b$ and thus
commutes with $d$ and $e$.  It remains to show that it commutes with
$c$.  But
\begin{align*}
\mathtt{c[[aa,BB],BB]C}&\mathtt{=[[aa,BB],cBBC]}\\
\mathtt{c[[aa,BB],BB]C}&\mathtt{=[[aa,BB],BB[c,BB]]}
\end{align*}

and $[c,b]$ commutes with $[a,b]$.  Thus $[d,[b,c]]$ is central, as is
$[b,[c,d]]$ and thus all 4-element commutators can be filled with a
horizontal disc.
\end{itemize}

Finally, by the Jacobi identity,
$$[d,[b,c]]=g_1^{-1}[b,[c,d]]^{-1}g_1g_2^{-1}[c,[b,d]]g_2=[b,[c,d]]^{-1},$$ so there is only one non-trivial 3-element
commutator.  Thus, any word representing the identity in $H$ can be
filled with a horizontal disc.  In particular, given a presentation of
$H$ with generators $a,b,c,d,e$, the relators can be filled with
horizontal discs; this gives a quadratic bound on the area of scalings
of the relators and thus, by Theorem \ref{mainscale}, a quadratic
isoperimetric inequality.

\end{document}